\documentclass[12pt,reqno]{amsart}
\usepackage[margin=1.1in]{geometry}
\usepackage{amsmath, amssymb, amsfonts}
\usepackage{mathtools}
\usepackage{hyperref}
\usepackage{ mathrsfs}
\usepackage{xcolor}
\def\sH{\mathscr{H}}
\def\sA{\mathscr{A}}
\def\sI{\mathscr{I}}
\def\sE{\mathscr{E}}

\def\eps{\varepsilon}
\def\ex{\mathrm{ex}}

\def\oF{\overrightarrow{F}}

\def\K4{K_4^{(3)}}
\def\Kr{K_r^{(r-1)}}
\def\K*{K^{(r-1)}_{a_1, \dots, a_r}}

\newtheorem{theorem}{Theorem}

\newtheorem{proposition}[theorem]{Proposition}

\newtheorem{lemma}[theorem]{Lemma}
\newtheorem{claim}[theorem]{Claim}

\newtheorem{fact}[theorem]{Fact}
\newtheorem{question}[theorem]{Question}

\begin{document}

\title{The number of cliques in hypergraphs with forbidden subgraphs}
\thanks{The first and second authors are partially supported by NSF grants DMS 1764385 and 2300347.
The third author is partially supported by NSF grant DMS  2300346 and Simons Collaboration Grant 710094.}
\author{Ayush Basu}
\author{Vojt\v{e}ch R\"odl}
\address{Department of Mathematics, Emory University, 
    Atlanta, GA, USA}
\email{\{abasu|vrodl\}@emory.edu}

\author{Yi Zhao}
\address{Department of Mathematics and Statistics, Georgia State University, Atlanta, GA 30303}
\email{yzhao6@gsu.edu}
\date{\today}
\subjclass{05C65, 05C35}
\keywords{Generalized Tur\'an problem, hypergraph, removal lemma}

\begin{abstract}
We study the maximum number of $r$-vertex cliques in $(r-1)$-uniform hypergraphs not containing complete $r$-partite hypergraphs $K_r^{(r-1)}(a_1, \dots, a_r)$. By using the hypergraph removal lemma, we show that this maximum is $o( n^{r - 1/(a_1 \cdots a_{r-1})} )$. This immediately implies the corresponding results of Mubayi and Mukherjee and of Balogh, Jiang, and Luo for graphs. We also provide a lower bound by using hypergraph Tur\'an numbers.
\end{abstract} 

\maketitle

\section{Introduction}
Given integers $r\geq 2$ and $n>0$ and two $r$-uniform hypergraphs $T$ and $F$, let $\ex(n, T, F)$ denote the maximum number of copies of $T$ in any $F$-free (i.e., not containing $F$ as a subgraph) $r$-uniform hypergraph on $n$ vertices. The case when $T$ is an edge (i.e., $T= K_2$) is the Tur\'an problem $\ex(n, F)$. The parameter $\ex(n,T,F)$ has been studied for different choices of graphs $T$ and $F$ by many authors (for example, see \cite{MR3548290, balogh2024maximum, MR0151956, MR3904833, MR4103873,MR3826292, MR645860, MR4546048}). 

Given an $r$-uniform hypergraph $F$ with vertex set $V(F)= \{v_1,\dots, v_\ell\}$, let $F (a_1,\dots, a_\ell)$ denote a \emph{blowup of $F$}, i.e., the hypergraph obtained from $F$ by replacing each vertex $v_i$ by a set $V_i$ of size $a_i$, and every edge  $\{v_{i_1},\dots, v_{i_r}\}$ by a complete $r$-partite $r$-uniform hypergraph on the vertex sets $V_{i_1},\dots, V_{i_r}$. If $a_1 = \cdots = a_\ell = a$, then we denote $F (a_1,\dots, a_\ell)$ by $F(a)$. 
For $r\le \ell$, we denote by $K_{\ell}^{(r)}(a_1, \dots, a_{\ell})$ the complete $\ell$-partite $r$-uniform hypergraph with $a_1, \dots, a_{\ell}$ vertices in its parts. We often omit the superscript when $r=2$, for example, $K_3$ is a graph triangle and $K_3(1, a, b)$ is a complete tripartite graph with parts of size $1$, $a$ and $b$.

In this note, we consider the parameter $\ex(n,T,F)$, when $T$ is an $r$-uniform hypergraph, and $F$ is a blowup of $T$. 
This problem is related to the following classical result of Erd\H os \cite{MR183654} on the Tur\'an number $\ex(n, K^{(r)}_r(a_1, \dots, a_r))$. 
It states that given integers $r\ge 2$ and $1\le a_1\le \cdots\le a_r$, 
\begin{align}
\label{eqn: Erdos}
\ex(n, K^{(r)}_r(a_1, \dots, a_r)) = O (n^{r - 1/(a_1 \cdots a_{r-1})} ).
\end{align}
Note that the $r=2$ case of \eqref{eq:Erdos} follows from a result of K\"ovari, S\'os, and Tur\'an~\cite{MR65617}.

Given $2\le s\leq r$ and an $r$-uniform hypergraph $F$, the \emph{$s$-uniform shadow} 
$\partial^{(s)} F$ of $F$ is an $s$-uniform hypergraph on $V(F)$ whose edge set consists all $s$-subsets $A\subseteq V(F)$ such that $A\subseteq B$ for some edge $B\in F$. We observe the following simple fact (and will prove it in Section \ref{sec2}).
\begin{fact}\label{fact1}
    Given $r\ge 3$ and an $r$-uniform hypergraph $F$, we have
    \begin{align*} 
\ex(n, K_r, \partial ^{(2)}F) \le \cdots \le \ex(n, K_r^{(r-1)}, \partial ^{(r-1)}F) \le \ex(n, F).
\end{align*}
\end{fact}
Fact~\ref{fact1} and (\ref{eqn: Erdos}) together imply that, given positive integers $r\ge 3$ and $a_1\le \cdots\le a_r$, 
\begin{align}
    \label{eq:Erdos}
    \ex(n, K^{(r-1)}_r, K_r^{(r-1)}(a_1,\dots , a_r)) \le \ex(n, K^{(r)}_r(a_1, \dots, a_r)) = O (n^{r - 1/(a_1 \cdots a_{r-1})} ).
\end{align}

Using the triangle removal lemma, Mubayi and Mukherjee \cite{MR4546048} showed that for any $1\le a\le b$,
$\ex(n, K_3, K_3(1, a, b)) = o(n^{3-1/a})$.\footnote{It was mentioned in \cite{balogh2024maximum, MR4546048} that this result was also proved by several other researchers.}
Our main result extends this to hypergraphs. 
\begin{theorem}
\label{prop:ub1}
    Given positive integers $r\geq 3$ and $ a_1\le \cdots \le a_r$, we have, 
    \begin{align}
    \label{eqn: theorempt1}
        \ex(n, K_r^{(r-1)}, K_r^{(r-1)}(a_1,\dots , a_r)) = o(n^{r - 1/(a_1 \cdots a_{r-1})}).
    \end{align}
    Further, given integers $a\geq 1$, $\ell \geq r$, and any $(r-1)$-uniform hypergraph $F$ on $\ell$ vertices, 
    \begin{align}
    \label{eqn: theorempt2}
        \ex(n, F, F(a)) = o\left(n^{\ell - \frac{1}{a^{\ell -1}}}\right).
    \end{align}
\end{theorem}
Very recently and independently, Balogh, Jiang, and Luo \cite{balogh2024maximum} proved that for any integers $r\geq 3$ and $1\le a_1\le \cdots \le a_r$, 
\begin{align}\label{eq:Balogh}
    \ex(n,K_r, K_r(a_1,\dots , a_r)) = o(n^{r - 1/(a_1 \cdots a_{r-1})}).
\end{align}
Note that Fact~\ref{fact1} and Theorem \ref{prop:ub1} together imply \eqref{eq:Balogh}.

Both (\ref{eqn: theorempt1}) and (\ref{eq:Balogh}) generalize the result of \cite{MR4546048} and all three proofs are similar in the sense that they all reduce the problem to an application of the removal lemma. While in \cite{balogh2024maximum}, the authors remarked that the proof idea in \cite{MR4546048} seemed not to work when proving (\ref{eq:Balogh}) for $r>3$ and $a_1 > 1$, here we present a proof of Theorem \ref{prop:ub1} (which also implies (\ref{eq:Balogh})) by taking an approach along the lines of \cite{MR4546048}.

The following lower bound complements Theorem \ref{prop:ub1} and will be proved in Section \ref{sec3}. Note that similar constructions (for graphs) appeared in \cite{balogh2024maximum, MR4546048}.
\begin{proposition}
\label{prop:lb4}
Given integers $1\le a_1\le \dots \le a_r$,
\begin{align}
\label{eqn: lb4}
    \ex(n, K_r^{(r-1)}, K_r^{(r-1)}(a_1,\dots , a_r)) =\Omega\left(n \cdot\ex\left( n, \, K^{(r-1)}_{r-1}(a_1, \dots, a_{r-1}) \right)\right).
\end{align}

\end{proposition}
Observe that if $a_1 = \cdots = a_{r-1}= 1$ and $a_r \geq 2$, then the right hand side of (\ref{eqn: lb4}) is zero and hence the lower bound is trivial. In this case, a construction in \cite{MR1884430} implies the following lower bound. Let $r_r(n)$ denotes the size of the largest subset of $[n]$ that does not contain an arithmetic progression of length $r$.
\begin{proposition}
\label{prop: lbAP}
    For every $r\geq 3$,
    \begin{align*}
    \ex(n, K_r^{(r-1)}, K_r^{(r-1)}(1,\dots,1,2)) \ge n^{r-2}r_r(n).
    \end{align*}
\end{proposition}
For the proof of Proposition~\ref{prop: lbAP}, see Section 3.


\section{Proof of Fact~\ref{fact1} and Theorem \ref{prop:ub1}}
\label{sec2}
In this section, we will prove Fact~\ref{fact1} and Theorem \ref{prop:ub1}. 

\begin{proof}[Proof of Fact~\ref{fact1}]
It suffices to show that, for every $2\le s\le r-1$,
\[
\ex(n, K_r^{(s)}, \partial^{(s)} F)\le \ex(n, K_r^{(s+1)}, \partial^{(s+1)} F),
\]
(trivially $\partial^{(r)} F=F$). Indeed, let $G$ be an $\partial^{(s)}F$-free $s$-graph on $[n]$ with $\ex(n, K_r^{(s)}, \partial^{(s)}F)$  copies of $K_r^{(s)}$. Let $H$ be the $(s+1)$-graph on $[n]$ whose edges are $(s+1)$-sets that span a copy of $K_{s+1}^{(s)}$ in $G$. 
We claim that $H$ is $\partial^{(s+1)}F$-free. Suppose instead, that $H$ contains 
a copy of $\partial^{(s+1)}F$ on some set $S\subset [n]$ under a bijection $\phi: V(F)\to S$. Consider an $s$-set $A\in \partial^{(s)}F$. We know $A\subset B$ for some $B\in \partial^{(s+1)}F$. Thus, $\phi(B)\in H$ by the definition of $\phi$ and consequently, $\phi(A)\in G$ by the definition of $H$. This implies that $S$ spans a copy of $\partial^{(s)}F$ in $G$, contradicting that $G$ is $\partial^{(s)}F$-free. 

Furthermore, it is easy to see that for any $r$-subset $S\subset [n]$, $S$ spans a copy of $K_r^{(s)}$ in $G$ if and only if $S$ spans a copy of $K_r^{(s+1)}$ in $H$. Thus, the number of $K_r^{(s+1)}$ in $H$ equals to $\ex(n, K_r^{(s)}, \partial^{(s)}F)$, the number of $K_r^{(s)}$ in $G$. Since $H$ is $\partial^{(s+1)}F$-free, we conclude that $\ex(n, K_r^{(s)}, \partial^{(s)}F)\le \ex(n, K_r^{(s+1)}, \partial^{(s+1)}F)$.
\end{proof}


Before proving Theorem \ref{prop:ub1}, we will fix some notation that we use for the rest of the section. We call $r$-uniform hypergraphs \emph{$r$-graphs}. Given an $(r-1)$-graph $G$ and a vertex $v\in V(G)$, let $G(v)$ be the $(r-1)$-graph with vertex set $V(G)\setminus \{v\}$, and
    \begin{align*}
      \{v_1,\dots, v_{r-1}\} \in G(v) \text{ if } \{v,v_1,\dots, v_{r-1}\} \text{ induces } K_r^{(r-1)} \text{ in } G.
    \end{align*}
    For a positive integer $a$, let $G(v_1)\cap \cdots \cap G(v_{a})$ be the $(r-1)$-graph with vertex set $V(G)\setminus \{v_1,\dots, v_{a}\}$ and edge set consisting of all $\{w_1,\dots, w_{r-1}\}$ such that $\{v_i, w_1,\dots, w_{r-1}\}$ induces a $K_{r}^{(r-1)}$ for every $i\in \{1,\dots, a\}$.
    
In the following proofs we will use the hypergraph removal lemma, which we state below.
\begin{lemma}[Hypergraph Removal Lemma \cite{MR2198495, MR2373376}]
    For every $r\geq 3$, $\varepsilon > 0$ there exists $\delta > 0$ such that for every $(r-1)$-uniform hypergraph $G$ on $n$ vertices the following holds. If $G$ contains at least $\varepsilon n^{r-1}$ edge disjoint copies of $K_r^{(r-1)}$, then it must contain at least $\delta n^r$ copies of $K_r^{(r-1)}$. 
\end{lemma}
We also need the following simple claim.
    \begin{claim}
    \label{ub2}
      For every positive integer $b$, $r\geq 3$ and  $(r-1)$-graph $G$ on $n$ vertices, the following holds. If $\sI$ is a collection of cliques $K_r^{(r-1)}$ of $G$ such that every edge $e\in G$ is contained in less than $b$ cliques of $\sI$, then $G$ contains at least $\frac{|\sI|}{r(b-1)}$ edge disjoint copies of $K_r^{(r-1)}$.
\end{claim}
\begin{proof}
    For $G$ and $\sI$ satisfying the above assumptions, let $\sI_1\subseteq \sI$ be a maximum collection of pairwise edge disjoint cliques $K_r^{(r-1)}$ in $\sI$ and let $\sE$ be the union of edge sets of the cliques in $\sI_1$. Clearly $|\sE| = r\cdot|\sI_1|$. Since by assumption, each edge $e\in \sE$ is contained in at most $b-1$ cliques $K_r^{(r-1)}$ in $\sI$, there are at most $(b-1)r|\sI_1|$ cliques in $\sI$ containing some edge of $\sE$. Due to the maximality of $\sI_1$, it follows that $|\sI|\le (b-1)r |\sI_1|$ and thus $|\sI_1|\geq \frac{|\sI|}{(b-1)r}$. 
   \end{proof}
Now we prove Theorem \ref{prop:ub1}. 
\begin{proof}[Proof of Theorem \ref{prop:ub1}]   
    Fix $r\geq 3$ and integers $a_1\le \cdots \le a_r$. We first consider the case when $a_{r-1} = 1$.
    Let $\eps > 0$, and let $G$ be a $K_{r}^{(r-1)}(1,\dots ,1, a_r)$-free $(r-1)$-graph on $n$ vertices, i.e., every edge of $G$ is in at most $(a_r - 1)$ copies of $K_r^{(r-1)}$. Assume by contradiction, that the collection $\sI$ of all $K_r^{(r-1)}$ in $G$ has size at least $\eps n^{r-1}$. In view of Claim \ref{ub2}, $G$ must contain at least $\eps' n^{r-1}$ edge disjoint copies of $K_r^{(r-1)}$ where $\eps' = ((a_r-1)r)^{-1}\eps$. By the hypergraph removal lemma, this implies that there exists some $\delta > 0$ (depending on $\eps')$ such that $G$ contains $\delta n^{r}$ copies of $K_r^{(r-1)}$. However, this contradicts (\ref{eq:Erdos}).
    
    Next, we consider the case where $a_{r-1}\geq 2$. Let $G$ be a $K_r^{(r-1)}(a_1,\dots, a_{r})$-free $(r-1)$-graph on $n$ vertices. First we will show that for every $\{v_1,\dots, v_{a_{r-1}}\}\subseteq V(G)$, the $(r-1)$-graph $G(v_1)\cap \cdots \cap G(v_{a_{r-1}})$ is $K_{r-1}^{(r-1)}(a_1,\dots, a_{r-2}, a_r)$-free. Indeed,
        given $\{v_1,\dots, v_{a_{r-1}}\}\subseteq V(G)$, assume by contradiction, that the $(r-1)$-graph $G(v_1)\cap \cdots \cap G(v_{a_{r-1}})$ contains a copy of $K_{r-1}^{(r-1)}(a_1,\dots, a_{r-2}, a_r)$ with the vertex set $V_1\sqcup V_2\dots\sqcup V_{r-2} \sqcup V_r$, where $|V_i| = a_i$. Let $V_{r-1} := \{v_1,\dots, v_{a_{r-1}}\}$. Then the $r$-partite graph on $V_1\sqcup V_2\cdots \sqcup V_r$ in $G$ forms a copy of $K_{r}^{(r-1)}(a_1,\dots, a_{r})$, a contradiction.

    Consequently, for every $\{v_1,\dots, v_{a_{r-1}}\}\subseteq V(G)$,
    \begin{align}
    \label{eqn: Neighborhood}
        |G(v_1)\cap \cdots \cap G(v_{a_{r-1}})|\leq \ex(n,K_{r-1}^{(r-1)}(a_1,\dots, a_{r-2}, a_r))= O\left(n^{r-1 - \frac{1}{a_1a_2\cdots a_{r-2}}}\right).
    \end{align}
    Our goal is using the above fact to obtain a large collection of edge disjoint $K_r^{(r-1)}$ in $G$. To this end we consider the family $\sA$, elements of which are collections of $a_{r-1}$ copies of $K_r^{(r-1)}$ that share an edge of $G$. More formally,
    \begin{align*}
        \sA &:= \{\{T_1,\dots, T_{a_{r-1}}\}: T_i \cong K_r^{(r-1)} \text{ and } T_1, T_2, \dots , T_{a_{r-1}} \text{ share an edge of } G\}.
    \end{align*}
    Next we give an upper bound on the size of $\sA$. Given \textit{any} element in $\sA$, there exists vertices $\{v_1,\dots, v_{a_{r-1}}\}\subseteq V(G)$, and an edge $e\in G$ (in particular, $e\in G(v_1)\cap \cdots \cap G(v_{a_{r-1}})$), such that $e\cup\{v_i\}$ form a $K_r^{(r-1)}$ for every $1\leq i\leq a_{r-1}$. Consequently, the cardinality of $\sA$ can be bounded by the number of pairs $(\{v_1,\dots, v_{a_{r-1}}\},e)$ with $e\in G(v_1)\cap \cdots \cap G(v_{a_{r-1}})$. Thus in view of (\ref{eqn: Neighborhood}), 
    \begin{align}
    \label{clm:edgesinH}
        |\sA| \leq {n\choose a_{r-1}}|G(v_1)\cap \cdots \cap G(v_{a_{r-1}})|
       &\leq n^{a_{r-1}}O\left(n^{r-1 - \frac{1}{a_1a_2\cdots a_{r-2}}}\right).
    \end{align}
   In order to prove (\ref{eqn: theorempt1}) of Theorem \ref{prop:ub1}, assume by contradiction, that $G$ contains $N = \Omega(n^{r - 1/(a_1 \cdots a_{r-1})})$ copies of $K_{r}^{(r-1)}$. We will find a collection $\sI$ of cliques $K_{r}^{(r-1)}$ in $G$ satisfying 
   \begin{align}
   \label{eqn: claimsI}
    |\sI| = \Omega(n^{r-1})\quad \text{ and } \quad \sI^{(a_{r-1})}\cap \sA = \emptyset,   
   \end{align} 
   i.e., for every $S\subseteq \sI$ with $|S| = a_{r-1}$, $S$ is not an element of $\sA$. Note that if $|\sA|\leq N/2$, then one can obtain $\sI$ from the collection of $K_r^{(r-1)}$ in $G$, by deleting a copy of $K_r^{(r-1)}$ for each element of $\sA$. Thus, $\sI^{(a_{r-1})}\cap \sA = \emptyset$ and $|\sI|$ is at least $N/2 = \Omega(n^{r - 1/(a_1 \cdots a_{r-1})})$ which is bigger than
   $\Omega(n^{r-1})$.  
   
    Now we consider the case where $|\sA|\geq N/2$. Let $\mathbf{I}$ be a random subset of copies of $K_{r}^{(r-1)}$ where each copy of $K_r^{(r-1)}$ in $G$ is chosen with probability $p>0$ independently. Let $\mathbf{I}^{(a_{r-1})}$ denote the collection of $a_{r-1}$-subsets of $\mathbf{I}$. We have that,
     \begin{align*}
         \mathbb{E}[|\mathbf{I}|- |\sA\cap \mathbf{I}^{(a_{r-1})}|] = pN- p^{a_{r-1}}|\sA|.
     \end{align*}
    Let $p$ be chosen such that, $p^{a_{r-1}}|\sA| = pN/2$, which implies
    \begin{align*}
       p = \left(\frac{N}{2|\sA|}\right)^{\frac{1}{a_{r-1}-1}}\leq 1, \text{ (since } N\leq 2|\sA|) \quad \text{ and } \quad pN = \frac{N^\frac{a_{r-1}}{a_{r-1}-1}}{(2|\sA|)^\frac{1}{a_{r-1}-1}}.
    \end{align*}
   Consequently, there exists a choice of $\sI'$ such that,
   \begin{align*}
       |\sI'|- |\sA\cap \sI'^{(a_{r-1})}|\geq \frac{pN}{2}.
   \end{align*} 
   Let $\sI\subseteq\sI'$ be the collection of $K_r^{(r-1)}$ of $G$ formed by deleting one $K_r^{(r-1)}$ in $\sI'$ from every $a_{r-1}$ subset in $\sA\cap \sI'^{(a_{r-1})}$. Consequently, $\sA\cap \sI^{(a_{r-1})}$ is empty. Further,
   \begin{align}
   \label{eqn: n^r-1}
       |\sI|\geq \frac{pN}{2}= \frac{1}{2}\frac{N^\frac{a_{r-1}}{a_{r-1}-1}}{(2|\sA|)^\frac{1}{a_{r-1}-1}}.
   \end{align}
   Using the value of $N$ (by assumption) and $|\sA|$ in (\ref{clm:edgesinH}), the exponent of $n$ in the RHS of (\ref{eqn: n^r-1}) is equal to
    \begin{align*}
       &\left( r - \frac{1}{a_1a_2a_3\cdots a_{r-1}}\right)\frac{a_{r-1}}{a_{r-1}-1} - \left(a_{r-1} + r - 1 -  \frac{1}{a_1\cdots a_{r-2}}\right)\frac{1}{a_{r-1}-1}\\
       &=\frac{a_{r-1} r -a_{r-1} - (r -1)}{a_{r-1}-1} = r-1,  
    \end{align*}
    which implies $  |\sI| = \Omega(n^{r-1}).$ Hence $\sI$ satisfies (\ref{eqn: claimsI}).
    
   Next we obtain a family of edge disjoint $K_r^{(r-1)}$ in $G$ from $\sI$. By construction, $\sI$ is a collection of cliques $K_r^{(r-1)}$ in $G$ such that every edge $e\in G$ is contained in less than $a_{r-1}$ cliques of $\sI$. In view of Claim \ref{ub2}, this implies that $G$ contains at least $|\sI|/r(a_{r-1}-1)$ edge disjoint copies of $K_r^{(r-1)}$. Since $|\sI| = \Omega(n^{r-1})$, this implies that $G$ contains $\Omega(n^{r-1})$ copies of edge disjoint $K_r^{(r-1)}$.  
   
   To summarise, this implies that given any $\eps > 0$, and $(r-1)$-graph $G$ on $n$ vertices that is $K_r^{(r-1)}(a_1,\dots, a_{r})$-free the following holds. Assuming by contradiction that $G$ contains $N = \eps n^{r - 1/(a_1 \cdots a_{r-1})}$ copies of $K_r^{(r-1)}$,  there exists some $\eps' > 0$ (depending only on $\eps, r, a_i$) such that $G$ contains $\eps' n^{r-1}$ edge disjoint copies of $K_r^{(r-1)}$. By the hypergraph removal lemma, this implies that there exists some $\delta>0$ (depending only on $\eps, r, a_i$) such that $G$ contains $\delta n^{r}$ copies of $K_r^{(r-1)}$. In view of (\ref{eq:Erdos}), however, this implies that $G$ contains $K_r^{(r-1)}(a_1,\dots, a_{r})$. Thus (\ref{eqn: theorempt1}) holds.

Now we prove the upper bound in (\ref{eqn: theorempt2}) on $\ex(n, F, F(a))$ for any given $(r-1)$-graph $F$. Label the vertices of $F$ $v_1,\dots, v_\ell$. Let $G$ be an $F(a)$-free $(r-1)$-graph on $n$ vertices, and assume by contradiction, that $G$ contains $N = \Omega(n^{\ell - \frac{1}{a^{\ell -1}}})$ copies of $F$.
    Given an $\ell$-partition of $V(G) = W_1\sqcup\cdots \sqcup W_\ell$, we call a set $X\subseteq V(G)$ \textit{crossing} if $|X\cap W_i|\leq 1$ for $1\leq i\leq \ell$.
    We call a copy of $F$ in $G$ on a vertex set $\{x_1,\dots, x_\ell\}$ \textit{aligned} with respect to $W_1\sqcup W_2\sqcup \cdots \sqcup W_\ell$ if
    \begin{enumerate}
        \item $x_i\in W_i$ for $i = 1,2,\dots, \ell$, and
        \item $x_i \mapsto v_i$ is an isomorphism. 
    \end{enumerate}
    We will denote such a copy by $\oF$. 
    A simple averaging argument yields that there exists a partition of $V(G) = W_1\sqcup\cdots \sqcup W_\ell$ with at least $\ell^{-\ell}N$ copies of $\oF$. 
    
     Let $\sH$ be an \textit{auxiliary} $\ell$-partite $(\ell-1)$-graph with vertex set $W_1\sqcup\cdots \sqcup W_\ell$. Let the edges of $\sH$ be those crossing $(\ell - 1)$-tuples that extend to a copy of $\oF$. Formally,
    \begin{align*}
        \sH = \bigsqcup_{i=1}^{\ell}\left\{(x_j)_{j\in [\ell]\setminus\{i\}}: \text{ there exists } x_i \in W_i \text{ such that } (x_1,\dots, x_\ell)\text{ is a copy of } \oF\right\}.
    \end{align*}
    Note that each \textit{aligned copy} $\oF$ in $G$ forms a $K_{\ell}^{(\ell -1)}$ in $\sH$. Consequently, the number of copies of $K_\ell^{(\ell-1)}$ in $\sH$ is at least $(\ell^{-\ell})N = \Omega(n^{\ell - \frac{1}{a^{\ell -1}}})$. 
    
    By the first part of Theorem \ref{prop:ub1}, this implies that $\sH$ contains a copy of $K_{\ell}^{(\ell -1)}(a)$ with vertex sets $U_i\subseteq W_i$ for $1\leq i\leq \ell$. Let $(x_1,\dots, x_\ell) \in U_1\times \cdots \times U_\ell$. Since $\ell \geq r$, for every edge $\{v_{i_1},\dots, v_{i_{r-1}}\}$ of $F$, there exists an $\ell-1$ subset $S\subseteq [\ell]$ such that $\{i_1,\dots, i_{r-1}\} \subseteq S$. By definition of $\sH$, the tuple $(x_s:s\in S)$ must extend to some copy of $\oF$, which implies $\{x_{i_1},\dots, x_{i_{r-1}}\}$ must be an edge in $G$. 
    
    Consequently, for every $(x_1,\dots, x_\ell) \in U_1\times \cdots \times U_\ell$, we have that the subgraph of $G$ induced by $\{x_1,\dots, x_\ell\}$ contains an aligned copy $\oF$. This implies that $G$ contains a copy of $F(a)$, contradicting the assumption that $G$ is $F(a)$-free. 
\end{proof}


\section{Lower Bound Constructions}
\label{sec3}
In this section, we will prove Propositions \ref{prop:lb4} and 
\ref{prop: lbAP}. 
\begin{proof}[Proof of Proposition \ref{prop:lb4}]
We construct an $(r-1)$-graph $H$ whose vertex set is partitioned into $A\sqcup B$ such that 
\begin{itemize}
\item $|A|= \lfloor n/r \rfloor$ and $|B|= \lceil (r-1)n/r \rceil$;
\item $H[B]$ is $K^{(r-1)}_{r-1}(a_1, \dots, a_{r-1})$-free and has $\ex( \lceil \frac{r-1}r n \rceil, K^{(r-1)}_{r-1}(a_1, \dots, a_{r-1}))$ edges;
\item every vertex of $A$ and every $(r-2)$-subset of $B$ form an edge and there are no other edges intersecting $A$. In other words, the link of every vertex in $A$ is the complete $(r-2)$-graph on the vertex set $B$.
\end{itemize}
The number of $\Kr$ is at least $\lfloor \frac nr \rfloor\ex( \lceil \frac{r-1}r n \rceil, K^{(r-1)}_{r-1}(a_1, \dots, a_{r-1}))$ because every vertex of $A$ together with any edge of $B$ forms a copy of $\Kr$. It remains to show that $H$ contains no $K_r^{(r-1)}(a_1,\dots , a_r)$. Assume by contradiction, it does. Since there is no edge containing two vertices from $A$, and $a_1\leq a_2\le \cdots \le a_r$, the subgraph induced by $H$ on $B$ needs to contain a $K^{(r-1)}_{r-1}(a_1, \dots, a_{r-1})$, thus contradicting the construction of $H$.  
\end{proof}
The proof of Proposition \ref{prop: lbAP} is based on a construction given in \cite{MR1884430}.
\begin{proof}[Proof of Proposition \ref{prop: lbAP}]
   In the proof of \cite[Proposition 2.1]{MR1884430}, it was shown that for every $r\geq 3$, there exists an $r$-partite $r$-graph $H$ with parts $V_1,\dots, V_r$ satisfying the following properties. 
\begin{enumerate}
    \item For every $\{x_1,\dots, x_{r-1}\}\subseteq V(H)$, there exists at most one edge in $H$ containing $\{x_1,\dots, x_{r-1}\}$. 
    \item For every collection of subsets $\{\{x_i,y_i\}\subseteq V_i: 1\le i\le r\}$, there exists $1\le i\le r$ such that $\{x_1,\dots, x_r\}\setminus\{x_i\}\cup \{y_i\}$ is not an edge of $H$. 
    \item $H$ has $(r-1)r n$ vertices and $n^{r-2}r_r(n)$ edges. 
\end{enumerate}
Let $G$ be the $(r-1)$-uniform shadow of $H$, i.e., $G=\partial^{(r-1)}H$. We claim that $G$ is $K_r^{(r-1)}(1,\dots,1,2)$-free and contains $n^{r-2}r_r(n)$ copies of $K_r^{(r-1)}$.  Since $G$ is the shadow of $H$, the number of copies of $K_r^{(r-1)}$ in $G$ is at least the number of edges in $H$.   

While the edges of $H$ correspond to a collection of edge disjoint cliques (``real cliques") in $G$, we will now show that $G$ contains no other cliques $K_r^{(r-1)}$. Assume by contradiction that $\{x_1,\dots, x_{r}\}$ induces such a ``fake clique" $K_r^{(r-1)}$, i.e., $\{x_1,\dots, x_{r}\}\notin H$ but induces a $K_r^{(r-1)}$ in $G$. Since every edge of this clique belongs to some ``real clique", for every $1\le i\le r$, there must exist $y_i\neq x_i$ in $V_i$ such that $\{x_1,\dots, x_r\}\setminus\{x_i\}\cup \{y_i\} \in H$, contradicting (2). Consequently, by (1), no two $K_r^{(r-1)}$ in $G$ share an edge and hence $G$ is $K_r^{(r-1)}(1,\dots,1,2)$-free.
\end{proof}

\section{Concluding remarks}
As mentioned earlier, in the case where $a_1=\cdots = a_{r-1}= 1$ and $a_r\geq 2$ the lower bound in (\ref{eqn: lb4}) is trivial. We ask if there are other sequences of integers $a_1\leq \cdots \leq a_r$ for which (\ref{eqn: lb4}) can be improved.   
\begin{question}
\label{question}
     Given integer $r\geq 3$, for what sequence of integers $ 1\le a_1\le \cdots \le a_r$,
    $$\ex(n, K_r^{(r-1)}, K_r^{(r-1)}(a_1,\dots , a_r))  \geq
n^{1+\eps}\cdot\ex\left( n, \, K^{(r-1)}_{r-1}(a_1, \dots, a_{r-1}) \right) $$
for some $\eps = \eps(n)> 0$?
\end{question}

The order of magnitude for $\ex(n, K_r^{(r-1)}, K_r^{(r-1)}(a_1,\dots , a_r))$ is not known in any non-trivial case. The case when $r\geq 3$ and $a_1 = \cdots = a_r=2$ is related to a problem of Erd\H{o}s, see, e.g., \cite{MR4328787,MR1723801,MR1928871}. Theorem \ref{prop:ub1} and Proposition \ref{prop:lb4}, together with the lower bound in \cite{MR4328787} imply that 
\begin{align*}
    \Omega\left(n^{r-\left\lceil\frac{2^{r-1}-1}{r-1}\right\rceil^{-1}} \right)\le \ex(n, K_r^{(r-1)}, K_r^{(r-1)}(2,\dots , 2)) \le o\left(n^{r - \frac{1}{2^{r-1}}}\right).
\end{align*}
It was conjectured in \cite{MR1909504},  that $\ex(n, K^{(r-1)}_{r-1}(a_1, \dots, a_{r-1}) )= \Omega (n^{r-1 - 1/(a_1\cdots a_{r-2})})$. This was confirmed for some cases in \cite{MR3718079,MR1909504,PohoataZakharov}. If this conjecture is true, then Theorem \ref{prop:ub1} and Proposition~\ref{prop:lb4} would imply that,  
\begin{align*}
 \Omega (n^{r - 1/(a_1\cdots a_{r-2})}) \leq   \ex(n, K_r^{(r-1)}, K_r^{(r-1)}(a_1,\dots , a_r)) \leq o(n^{r - 1/(a_1 \cdots a_{r-1})}).
\end{align*}
When $a_1 = \cdots = a_r = a\ge 2$, one can obtain that $\ex(n, K^{(r-1)}_{r-1}(a)) = \Omega\left(n^{r-1 - \frac{(r-1)(a-1)}{a^{r-1} - 1}}\right)$ by using the {probabilistic deletion method} \cite{MR0382007}. Together with Proposition~\ref{prop:lb4} and Theorem \ref{prop:ub1}, this gives  
\[ \Omega \left( n^{r - \frac{(r-1)(a-1)}{a^{r-1} - 1} }\right) \leq \ex(n, K_r^{(r-1)}, K_r^{(r-1)}(a)) \leq o\left(n^{r - \frac{1}{a^{r-1}}}\right).
\]
When $a_1 = 1$ and $a_2 = \cdots = a_r = a\geq 2$, instead of Proposition~\ref{prop:lb4}, one can employ the {deletion method} directly to an random $(r-1)$-uniform hypergraph on $n$ vertices by removing copies of $K_r^{(r-1)}(1,a,\dots, a)$. Together with Theorem \ref{prop:ub1}, this implies that
\[
\Omega \left( n^{r - \frac{r(r-1)}{a^{r-2}}} \right)\leq \ex(n, K_r^{(r-1)}, K_r^{(r-1)}(1, a, \dots, a)) \leq o(n^{r - \frac{1}{a^{r-2}}}).
\]
It would be interesting to improve the gaps in any of the above cases.

\section*{Acknowledgment}
The authors thank Jie Han and Sean Longbrake for valuable discussions. They also thank two anonymous referees for  helpful comments.


\begin{thebibliography}{10}

\bibitem{MR3548290}
N.~Alon and C.~Shikhelman.
\newblock Many {$T$} copies in {$H$}-free graphs.
\newblock {\em J. Combin. Theory Ser. B}, 121:146--172, 2016.

\bibitem{balogh2024maximum}
J.~Balogh, S.~Jiang, and H.~Luo.
\newblock On the maximum number of $r$-cliques in graphs free of complete
  $r$-partite subgraphs, arXiv:2402.16818, 2024.

\bibitem{MR4328787}
D.~Conlon, C.~Pohoata, and D.~Zakharov.
\newblock Random multilinear maps and the {E}rd\H{o}s box problem.
\newblock {\em Discrete Anal.}, pages Paper No. 17, 8, 2021.

\bibitem{MR0151956}
P.~Erd\H{o}s.
\newblock On the number of complete subgraphs contained in certain graphs.
\newblock {\em Magyar Tud. Akad. Mat. Kutat\'{o} Int. K\"{o}zl.}, 7:459--464,
  1962.

\bibitem{MR183654}
P.~Erd\H{o}s.
\newblock On extremal problems of graphs and generalized graphs.
\newblock {\em Israel J. Math.}, 2:183--190, 1964.

\bibitem{MR0382007}
P.~Erd\H{o}s and J.~Spencer.
\newblock {\em Probabilistic methods in combinatorics}, volume Vol. 17 of {\em
  Probability and Mathematical Statistics}.
\newblock Academic Press [Harcourt Brace Jovanovich, Publishers], New
  York-London, 1974.

\bibitem{MR3904833}
B.~Ergemlidze, A.~Methuku, N.~Salia, and E.~Gy\H{o}ri.
\newblock A note on the maximum number of triangles in a {$C_5$}-free graph.
\newblock {\em J. Graph Theory}, 90(3):227--230, 2019.

\bibitem{MR1884430}
P.~Frankl and V.~R\"{o}dl.
\newblock Extremal problems on set systems.
\newblock {\em Random Structures Algorithms}, 20(2):131--164, 2002.

\bibitem{MR4103873}
D.~Gerbner, E.~Gy\H{o}ri, A.~Methuku, and M.~Vizer.
\newblock Generalized {T}ur\'{a}n problems for even cycles.
\newblock {\em J. Combin. Theory Ser. B}, 145:169--213, 2020.

\bibitem{MR3826292}
L.~Gishboliner and A.~Shapira.
\newblock A generalized {T}ur\'{a}n problem and its applications.
\newblock In {\em S{TOC}'18---{P}roceedings of the 50th {A}nnual {ACM} {SIGACT}
  {S}ymposium on {T}heory of {C}omputing}, pages 760--772. ACM, New York, 2018.

\bibitem{MR2373376}
W.~T. Gowers.
\newblock Hypergraph regularity and the multidimensional {S}zemer\'{e}di
  theorem.
\newblock {\em Ann. of Math. (2)}, 166(3):897--946, 2007.

\bibitem{MR1723801}
D.~S. Gunderson, V.~R\"{o}dl, and A.~Sidorenko.
\newblock Extremal problems for sets forming {B}oolean algebras and complete
  partite hypergraphs.
\newblock {\em J. Combin. Theory Ser. A}, 88(2):342--367, 1999.

\bibitem{MR1928871}
N.~H. Katz, E.~Krop, and M.~Maggioni.
\newblock Remarks on the box problem.
\newblock {\em Math. Res. Lett.}, 9(4):515--519, 2002.

\bibitem{MR645860}
J.~Koml\'{o}s, J.~Pintz, and E.~Szemer\'{e}di.
\newblock A lower bound for {H}eilbronn's problem.
\newblock {\em J. London Math. Soc. (2)}, 25(1):13--24, 1982.

\bibitem{MR65617}
T.~K\"ovari, and V. T.~S\'os, V. T. and P.~Tur\'an. 
\newblock On a problem of {K}. {Z}arankiewicz
\newblock {\em Colloq. Math.} 3, 50--57, 1954.

\bibitem{MR3718079}
J.~Ma, X.~Yuan, and M.~Zhang.
\newblock Some extremal results on complete degenerate hypergraphs.
\newblock {\em J. Combin. Theory Ser. A}, 154:598--609, 2018.

\bibitem{MR1909504}
D.~Mubayi.
\newblock Some exact results and new asymptotics for hypergraph {T}ur\'{a}n
  numbers.
\newblock {\em Combin. Probab. Comput.}, 11(3):299--309, 2002.

\bibitem{MR4546048}
D.~Mubayi and S.~Mukherjee.
\newblock Triangles in graphs without bipartite suspensions.
\newblock {\em Discrete Math.}, 346(6):Paper No. 113355, 19, 2023.

\bibitem{MR2198495}
B.~Nagle, V.~R\"{o}dl, and M.~Schacht.
\newblock The counting lemma for regular {$k$}-uniform hypergraphs.
\newblock {\em Random Structures Algorithms}, 28(2):113--179, 2006.

\bibitem{PohoataZakharov}
C.~Pohoata and D.~Zakharov.
\newblock Norm hypergraphs.
\newblock {\em Combinatorica}, accepted, arXiv preprint arXiv:2101.00715.
\end{thebibliography}
\end{document}